\documentclass[12pt, draft]{amsart}
\usepackage[cp850]{inputenc}
\usepackage{amssymb}
\usepackage[mathscr]{eucal}
\usepackage{amscd}
\usepackage{amsmath}
\usepackage[all]{xy}

\newcommand{\R}{\mathbb R}
\newcommand{\N}{\mathbb N}

\newcommand{\F}{\mathscr{F}}

\newcommand{\KP}{{\sf K}\hspace{-1pt}{\sf P}}

\newcommand{\To}{\longrightarrow}

\newcommand{\Ran}{\mathop{\mathrm{Ran}}}

\newcommand{\Dom}{\mathrm{Dom}\,}

\theoremstyle{plain}
\newtheorem{theorem}{Theorem}[section]
\newtheorem{prop}[theorem]{Proposition}
\newtheorem{lemma}[theorem]{Lemma}
\newtheorem{cor}[theorem]{Corollary}

\newtheorem{defn}[theorem]{Definition}

\newcommand{\adef}{\begin{defn}}
	\newcommand{\zdef}{\end{defn}}

\title{The structure of Rochberg spaces}

\author[J.M.F. Castillo]{Jes\'us M. F. Castillo}
\address{Universidad de Extremadura, Instituto de Matem\'aticas Imuex,
E-06011 Badajoz, Spain.}
\email{castillo@unex.es}

\author[M. Gonz\'alez]{Manuel Gonz\'alez}
\address{Departamento de Matem\'aticas,
Universidad de Cantabria, E-39071 Santander, Spain.}
\email{manuel.gonzalez@unican.es}

\author[Ra\'ul Pino]{Ra\'ul Pino}
\address{Departamento de Matem\'aticas,
Universidad de Extremadura, E-06011 Badajoz, Spain.}
\email{rpino@unex.es}

\thanks{The research of the first and second authors has been supported  by MINCIN Project PID2019-103961. The research of the first and third authors has been supported  by Project IB20038 de la Junta de Extremadura.}


\begin{document}
\maketitle
\begin{abstract} We study the structure of the Rochberg Banach spaces $\mathfrak Z_n$ associated to the interpolation pair $(\ell_\infty, \ell_1)$ at $1/2$, and the operators defined on them.\end{abstract}

\section{Introduction}

Operators on the remarkable Kalton-Peck space $Z_2$ constructed in \cite{kaltpeck} have been studied in \cite{cgp}. But it turns out that in the same way that the space $Z_2$ is a nontrivial generalization of the Hilbert space $\ell_2$, the higher order Rochberg spaces $\mathfrak Z_n$ are nontrivial generalizations of both: $\mathfrak Z_1=\ell_2$ and $\mathfrak Z_2=Z_2$. This paper is devoted to study the Banach spaces $\mathfrak Z_n$ and  the operators defined on them. Some block operators on $\mathfrak Z_n$ have been obtained in \cite{gspaces} and \cite{symplectic}, while the peculiar structure of $\mathfrak Z_3$ was considered in \cite{sym}. Among other results we will show:

\begin{enumerate}
\item $\mathfrak Z_n$ is isomorphic to its dual (Proposition \ref{dualRochberg}). This result was obtained in \cite{ccc} but we will provide here an interpolation-free proof of independent interest.
\item Every normalized basic sequence in $\mathfrak Z_n$ admits a subsequence equivalent to one of the following $n$ types: the canonical basis of the Orlicz space $\ell_{f_k}$ generated by the function $f_k(t)= t^2 \log^{2(k-1)} t$, for $k=1,...,n$ (Theorem \ref{T_subsequence}). In particular, $\mathfrak Z_n$ admits exactly $n$ non-equivalent types of symmetric basic sequences.
\item $\mathfrak Z_n$ is not a complemented subspace of either a Banach lattice or a Banach space with GL.l.u.st; and does not contain complemented subspaces that are either Banach lattices or spaces with GL.l.u.st (Theorem \ref{lattice}).
\item Every operator $S:\mathfrak Z_n \To X$ is either strictly singular or an isomorphism on a subspace $E\subset \mathfrak Z_n$ isomorphic to
$\mathfrak Z_n$ and complemented in $\mathfrak Z_n$. Every operator $T: \mathfrak Z_n \To \mathfrak Z_n$ is either strictly singular or an isomorphism on a complemented copy $E$ of $\mathfrak Z_n$ such that $T[E]$ is complemented in $\mathfrak Z_n$ (Proposition \ref{seven}).
\item The space $\mathfrak L(\mathfrak Z_n)$ of operators on $\mathfrak Z_n$ admits a unique maximal ideal: that of strictly singular operators. Strictly singular and strictly cosingular operators of $\mathfrak L(\mathfrak Z_n)$ coincide (Theorem \ref{singcosingr}).
\item The composition of $n$ strictly singular operators of $\mathfrak L(\mathfrak Z_n)$ is compact, while the composition of $n-1$ operators is not necessarily compact (Theorem \ref{composition}).
\item Every copy of $\mathfrak Z_n$ in $\mathfrak Z_n$ is complemented (Theorem \ref{T_complementado}).
\item $\mathfrak Z_n$ is $\mathfrak Z_n$-automorphic; namely, every isomorphism between two infinite codimensional subspaces $A,B$ of $\mathfrak Z_n$ isomorphic to $\mathfrak Z_n$ can be extended to an automorphism of $\mathfrak Z_n$ (Proposition \ref{automorphic}).
\item A generalized version of \cite[Corollary 5.9]{ccfm}: if $T$ and $S$ are upper triangular operators on $\mathfrak Z_n$ making a commutative diagram
$$\xymatrix{
0 \ar[r] & \mathfrak Z_n \ar[r]\ar[d]_{T} & \mathfrak Z_{2n} \ar[r] \ar[d] & \mathfrak Z_n \ar[r]
\ar[d]^{S} & 0\\
0 \ar[r] & \mathfrak Z_n  \ar[r] & \mathfrak Z_{2n} \ar[r] & \mathfrak Z_n \ar[r] & 0}$$
then $T-S$ is strictly singular (Proposition \ref{compact}).\end{enumerate}

The final Appendix of the paper is devoted to  obtain a generalized commutator theorem valid for operators on arbitrary Rochberg spaces
(not only those generated from the scale of $\ell_p$-spaces) that extends results in \cite{caceso2} and \cite{rochberg}.

\section{Preliminaries on general Rochberg spaces}
We will work with complex interpolation of pairs $(X_0, X_1)$ of Banach spaces in the classical setting \cite{BL} that we briefly expose next:
The space of the pair $(X_0, X_1)$ are linear and continuously embedded in a larger Banach space $\Sigma$, which can be
assumed  to be $\Sigma = X_0+X_1$ endowed with the norm $\|x\|= \inf \{\|x_0\|_0 + \|x_1\|_1: x= x_0 + x_1\; x_j\in X_j \; \mathrm{for}\; j=0,1\}$. The pair will be called {\em regular} if, additionally, the intersection space $X_0\cap X_1$ is dense in both $X_0$
and $X_1$. We denote by $\mathbb S$ the complex strip defined by $0<Re(z)<1$. We will consider the Calder\'on space $\mathcal C(X_0,X_1)$ of continuous bounded functions $f:\overline{ \mathbb S }\To \Sigma$ that are holomorphic on $\mathbb S$ and satisfy the boundary condition that for $k=0,1$, $f(k+it)\in X_k$ for each $t\in\R$ and $\sup_t\|f(k+it)\|_{X_k}<\infty.$ The space  $\mathcal C(X_0,X_1)$ is complete under the norm
$\|f\| = \sup\{\|f(k+it)\|_{X_k}: k=0,1; t\in\R \}$. The family of interpolation spaces generated with this setting is
$$X_z=\{x\in {\Sigma}: x = f(z) \text{ for some } f\in\mathcal C(X_0, X_1)\}$$
endowed with its natural quotient norm.



Starting with $\mathcal C = \mathcal C(X_0, X_1)$, $n\in \N$ and $z \in \mathbb S$, the Rochberg space $\mathfrak R_n$ \cite{rochberg} is formed by the arrays of the truncated sequence of the Taylor coefficients of the elements of $\F$, namely $$\mathfrak R_n = \left\{\left (\frac{f^{(n-1)}(z)}{(n-1)!}, \dots, f^{(1)}(z), f(z)\right): f\in \mathcal C \right\}$$
endowed with the natural quotient norm. The space $\mathfrak R_1$ of arrays of length one (the values of the functions of $\mathcal C$ at $z$) correspond, in the suitable context, to classical interpolation spaces. The space $\mathfrak R_2$ of arrays of length two at $z$ constitute the usually called first derived space \cite{rochweiss}. It is clear that if $\delta_z:\mathcal C\to \Sigma$ is the canonical evaluation map then $\mathcal C/\ker{\delta_z}$ is isometric to $\mathfrak R_1 = \{w\in \Sigma: w=f(z) \text{ for some } f\in\mathcal C\}$, endowed with the quotient norm $\|w\|_{X_z}= \inf_{w=f(z)}\|f\|_\mathcal C$. If $\delta_z^{(n)} :\mathcal C\to \Sigma$ denotes the evaluation of the $n$-{th} derivative at $z$ then this map is bounded for all $z\in \mathbb S$ and all $n\in \N$ by the boundedness of $\delta_z$, the definition of derivative and the Banach-Steinhaus theorem (see \cite[Lemma 3.5]{jussieu} for an explicit estimation of its norm). We omit from now on the subscript $z$ and will consider it fixed. Consider the map $\Delta_n = \frac{1}{n!}\delta^{(n)}$. One has
$$\mathfrak R_n = \mathcal C/ \cap_{j<n}\ker\Delta_j = \left\{(w_{n-1},\dots,w_{0})\in \Sigma^n: w_i= \Delta_i f \text{ for some } f\in \mathcal C \text{ and all } 0\leq i<n \right\}.$$

The spaces $\mathfrak R_n$ can be arranged into exact sequences in a very natural way: this is implicit in \cite{rochberg} and explicit in \cite{cck}. Indeed, if for $1\leq n,
k<m$ we denote by $\imath_{n,m}:\Sigma^n\to \Sigma^m$ the inclusion on the
left given by
$\imath_{n,m}(x_n,\dots,x_1)=(x_n,\dots,x_1,0\dots,0)$ and by
$\pi_{m,k}:\Sigma^m\to \Sigma^k$ the projection on the right given by
$\pi_{m,k}(x_m,\dots, x_k,\dots,x_1)=  (x_k,\dots,x_1)$, then
$\pi_{m,k}$ restricts to an isometric quotient map of
$\mathfrak R_m$ onto $
\mathfrak R_k $ (this is trivial) and $\imath_{n,m}$ is an isomorphic embedding of
$\mathfrak R_n$ into $\mathfrak R_m$ (this can be proved as \cite[Proposition 2(a)]{cck}) and thus, see \cite[Theorem 4]{cck}, for each $n,k\in\mathbb N$ there is an exact sequence $$
\xymatrix{0\ar[r]& \mathfrak R_n \ar[r]^{\imath_{n,n+k}}& \mathfrak R_{n+k}\ar[r]^{\pi_{n+k,k}}&\mathfrak R_k\ar[r]& 0}.
$$
%

The sequences can be described via quasilinear maps $\Omega_{k,n}: \mathfrak R_k \To \Sigma^n$. To define them, let us set the map
$ \tau_{(n,0]}(f)=\left( \Delta_{n-1} f,\dots, \Delta_0 f \right).$ Fix $\varepsilon\in(0,1)$ and, for each $x=(x_{k-1},\dots,x_{0})$ in
$\mathfrak R_k$, select $f_x\in \mathcal C$ such that $x=\tau_{(k,0]}f_x(z)$, with $\|f_x\|\leq (1+\varepsilon)\|x\|$, in such a way that $f_x$ depends homogeneously on $x$. Then define
$\Omega_{k,n}(x)=\tau_{(n+k,k]}f_x(z)$.

It is clear that this map depends on the choice of $f_x$, but different choices of $f_x$ only produce bounded perturbations of the same map.  Any $\Omega_{k,n}$ defined in this way is a quasilinear map from $\mathfrak R_k$ to $\mathfrak R_n$, which means that there is a constant $C$ such that, for every $x,y\in \mathfrak R_k$ the difference
$
\Omega_{k,n}(x+y)-
\Omega_{k,n}(x)-
\Omega_{k,n}(y)
$ falls into $\mathfrak R_n$ and obeys an estimate
$$
\|\Omega_{k,n}(x+y)-
\Omega_{k,n}(x)-
\Omega_{k,n}(y)\|_{\mathfrak R_n}\leq C\left(\|x\|_{\mathfrak R_k} + \|y\|_{\mathfrak R_k}	\right).
$$
The map $\Omega_{k,n}$ can be used to form the twisted sum space $$
\mathfrak R_n \oplus_{\Omega_{k,n}} \mathfrak R_k = \left\{(y,x)\in \Sigma^{n+k}: y-\Omega_{k,n}(x)\in \mathfrak R_n, x\in \mathfrak R_k\right\},$$ endowed with the quasinorm $\|(y,x)\|_{\Omega_{k,n}}= \left\| y-\Omega_{k,n}(x)\right\|_{\mathfrak R_n}+ \| x\|_{\mathfrak R_k}.$
It turns out that $\mathfrak R_n \oplus_{\Omega_{k,n}} \mathfrak R_k$ and $\mathfrak R_{n+k}$ are the same space, and that the quasinorm above is equivalent to the norm of $\mathfrak R_{n+k}$. Rochberg spaces have been studied in \cite{cck,ccc}.\\

\section{The Rochberg spaces $\mathfrak Z_n$ associated to the scale of $\ell_p$-spaces}

The Rochberg spaces associated to the scale of $\ell_p$ spaces at  $z=1/2$ are specially interesting: we can think of higher order Rochberg spaces as a natural generalization of $\ell_2$ and of the Kalton-Peck space $Z_2$. For this reason we will call
$\mathfrak Z_n$ the $n^{th}$-Rochberg space, so that $\mathfrak Z_1=\ell_2$ and $\mathfrak Z_2=Z_2$. Following the notation in \cite{hmbst}, we will call $\KP_{k,m}$ the associated differential (former $\Omega_{k,m}$). It has been calculated in \cite{ccc} that
$$\KP_{1, n}(x)= x\left(\frac{2^{n}}{n!}\log^{n}\left(\frac{|x|}{\|x\|_2}\right),
\dots,  2\log\left(\frac{|x|}{\|x\|_2}\right) \right).\medskip$$

\subsection{Duality}

The following result was proved in \cite{ccc} using interpolation techniques. Let us provide an interpolation-free proof of independent interest. The reader is warned that for this proof the elements are labeled $(x_1, \dots, x_n)$ instead of $(x_{n-1}, \dots, x_0)$ for logistic reasons (the duality is easier to describe that way).

\begin{prop}\label{dualRochberg} For each $n\geq1$ there is an isomorphism $D_n:\mathfrak Z_n\rightarrow\mathfrak Z_n^*$ given by
	$$D_n(x_{1},\ldots,x_n)(y_{1},\ldots,y_n) = \sum_{i+j=n+1}(-1)^{i}\langle x_i,y_j\rangle.$$
\end{prop}
\begin{proof} Let $\mathfrak Z_n^0$ denote the subspace of $\mathfrak Z_n$ formed by all the finitely supported elements and observe that
the map $D_n$ above is well defined as a map $\mathfrak Z_n^0\rightarrow\mathfrak Z_n^*$. The plan of the proof is to show that this restriction is injective, bounded and has dense range; then, that its extension to the whole $\mathfrak Z_n$ has closed range, which makes it onto.\medskip

Consider the bilinear map $\langle \cdot, \cdot \rangle :\mathfrak Z_n^0\times\mathfrak Z_n^*\rightarrow\mathbb{C}$ given by
	$$\langle (x_{1},\ldots,x_n),(y_{1},\ldots,y_n)\rangle =\sum_{i+j=n+1}(-1)^{i}\langle x_i,y_j\rangle,$$
and define for each $(x_1,\ldots,x_n)\in\mathfrak Z_n^0$ the following norm
$$\omega_n(x_1,\ldots,x_n)=\sup_{\|(x_1',\ldots,x_n')\|_{\mathfrak Z_n}\leq1}\Big|D_n(x_1,\ldots,x_n)(x_1',\ldots,x_n')\Big|.$$
The norm $\omega_n$ verifies the following properties:
\begin{enumerate}
\item $\omega_n(x,0,0,\ldots,0)=\|x\|_{\ell_2}$ for any $x\in\ell_2$.

\begin{proof} Just by considering the vectors $(\KP_{1,n-1}x', x')$ we obtain
		$$\omega_n(x,0,\ldots,0)=\sup_{\|(\KP_{1,n-1}x',x')\|_{\mathfrak Z_n}\leq 1} |\langle x,x'\rangle| =\sup_{\|x'\|_2\leq 1}|\langle x,x' \rangle|=\|x\|_2.\qedhere$$
	\end{proof}

\item $\omega_n(x_1,\ldots,x_{n-1},0)=\omega_{n-1}(x_1,\ldots,x_{n-1})$.
	\begin{proof}\begin{eqnarray*}
		&\;&\omega_n(x_1,\ldots,x_{n-1},0)\\&=&\sup_{\|(x_1',\ldots,x_n')\|_{\mathfrak Z_n}\leq1}\Big\{\Big|\langle x_1,x_n' \rangle -\langle x_2,x_{n-1}' \rangle+\cdots+(-1)^{n}\langle x_{n-1},x_2' \rangle\Big| \Big\}\\
		&=&\sup_{\|(\KP_{n-1,1}(x_1',\ldots,x_n'),x_2',\ldots,x_n')\|_{\mathfrak Z_n}\leq1}\Big\{\Big| \langle x_1,x_n' \rangle -\langle x_2,x_{n-1}' \rangle+\cdots+(-1)^{n}\langle x_{n-1},x_2' \rangle\Big| \Big\}\\
		&=&\sup_{\|(x_2',\ldots,x_n')\|_{\mathfrak Z_{n-1}}\leq1}\Big\{\Big| \langle x_1,x_n' \rangle -\langle x_2,x_{n-1}' \rangle+\cdots+(-1)^{n}\langle x_{n-1},x_2' \rangle \Big| \Big\}\\
		&=&\omega_{n-1}(x_1,\ldots,x_{n-1}).\qedhere
		\end{eqnarray*}\end{proof}

\item $\|x_n\|_{\ell_2}\leq \omega_n(x_1,\ldots,x_n)$ for $(x_1,\ldots,x_n)\in\mathfrak Z_n^0$.
	\begin{proof}
		\begin{eqnarray*}
		\omega_n(x_1,\ldots,x_n)&=&\sup_{\|(x_1',\ldots,x_n')\|_{\mathfrak Z_n}\leq1}\Big\{\Big|\langle x_1,x_n' \rangle -\langle x_2,x_{n-1}' \rangle+\cdots+(-1)^{n}\langle x_{n-1},x_2' \rangle\Big| \Big\}\\
		&\geq& \sup_{\|(x_1',0,\ldots,0)\|_{\mathfrak Z_n}\leq1}\Big\{\Big|\langle x_1,x_n' \rangle \Big| \Big\}=\sup_{\|x_1'\|_{\ell_2}\leq 1}\Big\{\Big|\langle x_1,x_n' \rangle \Big| \Big\}=\|x_n\|_{\ell_2}.\qedhere
		\end{eqnarray*}
	\end{proof}

The next inequality appears essentially proved in \cite[Lemma 5.3]{ccc} in a more general context. See also \cite[Apendix I]{Benyamini_Lindenstrauss}. We will denote by $\KP^k x =\frac{2^k}{k!}x\log^k\frac{|x|}{\|x\|}$ the $k$-th component of the differential map $\KP_{1,n-1}$. \medskip

\item\label{L_bound} $\Big|\sum_{k=0}^{n-1} (-1)^k \langle \KP^{n-1-k}x,\KP^kx' \rangle   \Big| \leq 2^{n-1} \|x\|_{\ell_2}\,\|x'\|_{\ell_2}.$
\begin{proof}
Observe that the left side equals
$$(*)=\Big|\sum_{k=0}^{n-1}(-1)^k\big\langle \Delta_{n-1-k}(B_{1/2}(x)),\Delta_k(B_{1/2}(x')) \big\rangle  \Big|$$
where  $B_{1/2}(x)=x|x|^{2z-1}$. If $f(z)=\langle B_{1/2}(x),B_{-1/2}(x')\rangle=\langle x|x|^{2z-1},x'|x'|^{2(1-z)-1} \rangle$ then $f$ is well defined because, for each $z\in\mathbb{S}$, $B_{1/2}(x)(z)$ and $B_{-1/2}(x')(z)$ belong, respectively, to some $\ell_p$ space and its dual $\ell_{p^*}$. Thus,
\begin{align*}
(*)&=\Big|\sum_{k=0}^{n-1}(-1)^k\frac{1}{(n-1-k)!}\frac{1}{k!}\big\langle \delta_{1/2}^{n-1-k}B_{1/2}(x),\delta_{1/2}^kB_{1/2}(x')\big\rangle   \Big|\\
&=\frac{1}{(n-1)!}\Big|\sum_{k=0}^{n-1}(-1)^k{{n-1}\choose k}\big\langle \delta_{1/2}^{n-1-k}B_{1/2}(x),\delta_{1/2}^kB_{1/2}(x')\big\rangle     \Big|\\
&=\frac{1}{(n-1)!}\Big|f^{(n-1)}(z)|_{z=1/2}\Big|.
\end{align*}

The function $f:\mathbb{S} \rightarrow\mathbb{C}$ is analytic and thus, by Cauchy integral formula (see \cite[Lemma 3.8]{jussieu}),
\begin{align*}
\frac{1}{(n-1)!}|f^{(n-1)}(1/2)| &=\frac{1}{(n-1)!}\Big|\frac{(n-1)!}{2\pi i}\int_{\Gamma}\frac{f(w)}{(w-1/2)^n}dw\Big|\\
&\leq \frac{1}{2\pi}\int_\Gamma \frac{|f(w)|}{(|w-1/2|^n)}d|w|\\
&\leq 2^{n-1}\,\|B_{1/2}(x)\|_\mathcal{C}\,\|B_{-1/2}(x')\|_\mathcal{C}\\
&\leq 2^{n-1}\,\|x\|_{\ell_2}\,\|x'\|_{\ell_2}= 2^{n-1}\|x\|_{\ell_2}\,\|x'\|_{\ell_2}.
\end{align*}
The last two inequalities are a consequence of the definition of $B_{1/2}(x)(z)$ and the fact that $\|B_{1/2}(x)\|_\mathcal{C}= \|x\|_{\ell_2}$. \end{proof}

\item\label{L_equivalence} 	$\omega_n(\KP^{n-1}x,\ldots,\KP x,x)\sim \|x\|_{\ell_2}$.
	\begin{proof}
The $\geq$ inequality:
\begin{eqnarray*}
&\;&\omega_n(\KP^{n-1}x,\ldots,\KP x,x)\\&=&\sup_{\|(x_1',\ldots,x_n')\|_{\mathfrak Z_n}\leq1}\Big\{\Big|\langle \KP^{n-1}x,x_n' \rangle -\langle \KP^{n-2}x,x_{n-1}' \rangle+\cdots+(-1)^{n+1}\langle x,x_1' \rangle\Big| \Big\}\\
&\geq&  \sup_{\|(x_1',0,\ldots,0)\|_{\mathfrak Z_n}\leq1}\Big\{\Big|\langle \KP^{n-1}x,0 \rangle -\langle \KP^{n-2}x,0 \rangle+\cdots+(-1)^{n+1}\langle x,x_1' \rangle\Big| \Big\}\\
&=&\sup_{\|x_1'\|_2\leq 1}\Big\{\Big|\langle x,x_n' \rangle  \Big|  \Big\}=\|x\|_{\ell_2}.
\end{eqnarray*}
To obtain the $\leq$ inequality, observe that:

		\begin{eqnarray*}
&\;&\omega_n(\KP^{n-1}x,\ldots,\KP x,x)\\
&=&\sup_{\|(x_1',\ldots,x_n')\|_{\mathfrak Z_n}\leq1}\Big\{\Big|\langle \KP^{n-1}x,x_n' \rangle -\langle \KP^{n-2}x,x_{n-1}' \rangle+\cdots+(-1)^{n+1}\langle x,x_1' \rangle\Big| \Big\}\\
&=&\sup_{\|(x_1',\ldots,x_n')\|_{\mathfrak Z_n}\leq1}\Big\{\Big|\langle \KP^{n-1}x,x_n' \rangle-\big[\langle \KP^{n-2}x,x_{n-1}'-\KP(x_n') \rangle+\langle \KP^{n-2}x,\KP x_n' \rangle\big]+\\
&+&\cdots+(-1)^{n+1}\big[\langle x,x_1'-\KP^{n-1}x_n' \rangle + \langle x,\KP^{n-1}x_n' \rangle \big]   \Big|   \Big\}\\
&\leq& \sup_{\|(x_1',\ldots,x_n')\|_{\mathfrak Z_n}\leq1}\Big\{\Big|\langle \KP^{n-1}x,x_n' \rangle-\langle \KP^{n-2}x,\KP(x_n') \rangle+ \cdots+(-1)^{n+1}\langle x,\KP^{n-1}x_n'\rangle \Big|\Big\}\\
&+&\sup_{\|(x_1',\ldots,x_n')\|_{\mathfrak Z_n}\leq1}\Big\{\Big|-\langle \KP^{n-2}x,x_{n-1}'-\KP(x_n') \rangle
+\cdots+(-1)^{n+1}\langle x,x_1'-\KP^{n-1}x_n' \rangle\Big|   \Big\}\\
&\leq& 2^{n-1}\|x\|_{\ell_2}+ 2^{n-2}\|x\|_{\ell_2},
\end{eqnarray*}
where in the last inequality we have used (\ref{L_bound}) twice. Indeed, the second  summand can be bounded from above by $2^{n-2}\|x\|_{\ell_2}$: the elements $y_k'=x_k'-\KP^{n-k}(x_n')$ form a vector $(y_1',\ldots,y_{n-1}')\in\mathfrak Z_{n-1}\subset \mathfrak Z_n$. We then deduce by induction on $n$ that
\begin{eqnarray*}&\;&\sup_{\|(y_1',\ldots, y_{n-1}')\|_{\mathfrak Z_{n-1}}\leq1} \Big\{\Big|- \langle \KP^{n-2}x,y_{n-1}' \rangle +\cdots+ (-1)^{n+1}\langle x,y_1' \rangle\Big| \Big\}\\
&=&\omega_{n-1}(\KP^{n-2}x,\ldots,\KP x,x) \sim  \|x\|_{\ell_2}.\qedhere\end{eqnarray*}
\end{proof}

\item $\omega_n$ is equivalent to $\|\cdot\|_{\mathfrak Z_n}$.

\begin{proof}
We will use induction on $n$. Kalton proved in \cite{kaltsym} the case $n=2$, so we may suppose that $\omega_{n-1}(x_1,\ldots,x_{n-1})\sim \|(x_1,\ldots,x_{n-1})\|_{\mathfrak Z_{n-1}}$. If $(x_1,\ldots,x_n)\in\mathfrak Z_n^0$ then
\begin{itemize}
\item $(x_1-\KP^{n-1}(x_n),\ldots,x_{n-1}-\KP(x_n),0)\in\mathfrak Z_{n-1}$, and
 \item $(\KP^{n-1}(x_n),\ldots,\KP(x_n),x_n)\in\mathfrak Z_n$.\end{itemize} Thus
	\begin{eqnarray*}
	&\;& \omega_n(x_1,\ldots,x_n)\\&=&\omega_n\Big((x_1,\ldots,x_{n-1},0)-(\KP^{n-1}(x_n),\ldots,\KP(x_n),0)\\
	&\quad&\quad \quad +(\KP^{n-1}(x_n),\ldots,\KP(x_n),0)+(0,\ldots,0,x_n)\Big)\\
	&\leq& \omega_n\big(x_1-\KP^{n-1}(x_n),\ldots,x_{n-1}-\KP(x_n),0\big)+\omega_n(\KP^{n-1}(x_n),\ldots,\KP(x_n),x_n)\\
&= & (**)\end{eqnarray*}
A combination of (*) and (\ref{L_equivalence}) yields
$$(**) \sim \omega_{n-1}(x_1-\KP^{n-1}(x_n),\ldots,x_{n-1}-\KP(x_n))+\|x_n\|_{\ell_2}$$
and, by induction hypothesis, we get
\begin{eqnarray*}
&\leq& C \|\big(x_1-\KP^{n-1}(x_n),\ldots,x_{n-1}-\KP(x_n)\big)\|_{\mathfrak Z_{n-1}}+\|x_n\|_{\ell_2}\\
	&\leq& C \|(x_1,\ldots,x_{n-1})-\KP_{1,n-1}(x_n)\|_{\mathfrak Z_{n-1}}+C\|x_n\|_{\mathfrak Z_{n-1}}\\
	&=&C\,\|(x_1,\ldots,x_n)\|_{\mathfrak Z_n}.
	\end{eqnarray*}
	
Conversely, since $\|(x_1,\ldots,x_n)\|_{\mathfrak Z_n} = \|(x_1,\ldots,x_{n-1})-\KP_{1,n-1}(x_n)\|_{\mathfrak Z_{n-1}}+\|x_n\|_{\ell_2}$, by induction hypothesis we get:
	$$\|(x_1,\ldots,x_n)\|_{\mathfrak Z_n}\sim C\, \omega_{n-1}\big((x_1,\ldots,x_{n-1})-\KP_{1,n-1}(x_n)\big)+\|x_n\|_{\ell_2}$$
which, using (*), is equal to $C\,\omega_n\big((x_1,\ldots,x_{n-1},0)-(\KP_{1,n-1}(x_n),0)\big)+\|x_n\|_{\ell_2}$ and then, by (3), to $$C\,\omega_n\Big((x_1,\ldots,x_{n-1},0)+(0,0\ldots,0,x_n)-(0,0\ldots,0,x_n)-(\KP_{1,n-1}(x_n),0)\Big)+\omega_n(x_1,\ldots,x_n).$$ This last term is bounded by $C\Big( \omega_n(x_1,\ldots,x_n)+\omega_n(\KP_{1,n-1}(x_n),x_n)\Big)+\omega_n(x_1,\ldots,x_n)$ and thus, applying (3) and (\ref{L_equivalence}), by $3C \omega_n(x_1,\ldots,x_n)$.\end{proof}
\end{enumerate}
We are ready to conclude the proof of Proposition \ref{dualRochberg}: the equivalence between $\omega_n$ and $\|\cdot\|_{\mathfrak Z_n}$
allows to define $D_n: \mathfrak Z_n \To \mathfrak Z_n^*$. This map is linear and injective. It has dense range because, otherwise, there is $\varphi\in \mathfrak Z_n^{**}=\mathfrak Z_n$ such that $D_n(x)(\varphi) =0$ for al $x\in \mathfrak Z_n$. Since $D_n(x)(\varphi)=D_n(\varphi)(x)$, it turns out that $D_n(\varphi)=0$, which implies $\varphi=0$. Finally, the equivalence $\|\cdot\|_{\mathfrak Z_n}\sim \omega_n(\cdot)$ entails that $D_n$ has closed range and, therefore, it must be surjective. Moreover, the previous estimates yield $\|D_n\|\leq 1$ and $\|D^{-1}_n\|\leq 3^{n-1}$.\end{proof}

It seems that the same proof should work for Rochberg spaces derived from an interpolation scale $(X,X^*)$ in which $X$ is a superreflexive K\"othe function space and $(X,X^*)_{1/2}=L_2$. One should change the duality and use the one described in \cite[Corollary 4]{cabcent}.

\section{Basic sequences in $\mathfrak Z_n$}
In this section we will work with the pair $(\ell_1, c_0)$ instead of $(\ell_1, \ell_\infty)$. It is well known \cite{BL}, see also \cite{ccc}, that the two pairs provide the same Rochberg spaces. The advantage of working with $(\ell_1, c_0)$ is that the sequence of vectors $(e_n)$ is a symmetric basis in both spaces, as well as in all interpolated spaces. More generally:

\begin{prop}\label{basis}
For all $1\leq k\leq n$ the canonical unit vectors form a symmetric basis for the Banach space
$$X_{n-k}= \Delta_k\left (\bigcap_{0\leq j\neq k \leq n-1} \ker\Delta_j \right).$$
\end{prop}
\begin{proof} Observe that for $k\neq 0$ one has
$$\Delta_k\left (\bigcap_{0\leq j\neq k \leq n-1} \ker\Delta_j \right) =
\Delta_{k-1}\left (\bigcap_{0\leq j\neq k-1 \leq n-1}  \ker\Delta_j \right).$$

Thus, an induction argument on $n$ plus the standard fact that $X_{0}=\ell_2$ in $\mathfrak R_n$ show that we only need to worry about $X_n =  \{x\in \mathfrak Z_n: \KP_{1,n-1} x \in \mathfrak Z_{n-1}\}$, called the \emph{Domain} of $\KP_{1,n}$ and denoted $\Dom \KP_{1,n-1}$. Both the natural projection operator $P_n$ onto the subspace generated by $\{e_1,\ldots, e_n\}$ and the operators induced by permutations of the basis are operators on the scale. Therefore, such operators are bounded on $\Dom \KP_{1,n-1}$ by \cite[Theorem 7.1]{racsam} (a proof is given in Proposition \ref{L_dominio} of the Appendix). Clearly $(e_n)$ is contained in $X_j$ and generates a dense subspace. Since for each $x\in \textrm{span}\{e_n: n\in\N\}$, $P_n x$ converges to $x$ in $X_j$, it does for each $x\in X_j$. Thus $(e_n)$ constitutes a symmetric basis for $X_j$.\end{proof}

Kalton and Peck showed in \cite{kaltpeck} that $X_1$ in $\mathfrak Z_2$ is the Orlicz space generated by the Orlicz function $f_1(t)=t^2\log^2t$. In \cite{sym} it is shown that $X_{2}$ in $\mathfrak Z_3$ is the Orlicz space generated by the Orlicz function $f_2=t^2\log^4t$. In the general case we have:

\begin{prop}\label{basisorlicz}
For all $0\leq j\leq n-1$ the space $X_j$ is isomorphic to the Orlicz space $\ell_{f_j}$ generated by the Orlicz function
$f_j(t) = t^2 \log^{2j}t$.
\end{prop}
\begin{proof} To ease the process, recall that, for normalized $x$, we have
$$\KP_{1, n}(x)= x\left(\frac{2^{n}}{n!}\log^{n}\left(\frac{|x|}{\|x\|_2}\right),
\dots,  2\log\left(\frac{|x|}{\|x\|_2}\right) \right) = x\left[\frac{2^j}{j!}\log^j |x|\right]_{j=n}^{1}$$
Observe that to estimate $\|(0,\dots, 0, x)\|_{\mathfrak Z_n}$ we need to perform $n-1$ steps until arriving to $\|\cdot\|_{\ell_2}$.

	Step 1. If we call $\bigstar = \|(0,\dots, 0, x)\|_{\mathfrak Z_n}$ then
	$$\bigstar = \|\KP_{1,n-1}(x)\|_{\mathfrak Z_{n-1}} = \left\| x\left[\frac{2^j}{j!}\log^j |x|\right]_{j=n-1}^{1}\right\|_{\mathfrak Z_{n-1}}$$
	
	Step 2.
	$$\left\| x\left[\frac{2^j}{j!}\log^j |x|\right]_{j=n-1}^{1}\right\|_{\mathfrak Z_{n-1}} = \left\| x\left[\frac{2^j}{j!}\log^j |x|\right]_{j=n-1}^{2}
	- \KP_{1,n-2}(\alpha_1 x \log |x|)\right\|_{\mathfrak Z_{n-2}} $$
	where $\alpha_1$ is the rightmost coefficient of $\left[\frac{2^j}{j!}\log^j |x|\right]_{j=n-1}^{1}$ (namely, 2). Thus, adding and subtracting
	$\alpha_1 \log |x| \KP_{1,n-2}(x)$ one gets
		
\begin{eqnarray*}
		\bigstar &\leq&
		\left\| x\left[\frac{2^j}{j!}\log^j |x|\right]_{j=n-1}^{2} - \alpha_1 \log |x| \KP_{1,n-2}(x)\right\|_{\mathfrak Z_{n-2}} \\
		&\quad\quad \quad +& \|\alpha_1 \log |x| \KP_{1,n-2}(x) -  \KP_{1,n-2}(\alpha_1 x \log |x|) \|_{\mathfrak Z_{n-2}}\\
		&\leq& \left\| x\left[\frac{2^j}{j!}\log^j |x|\right]_{j=n-1}^{2} - \alpha_1 \log |x| \KP_{1,n-2}(x)\right\|_{\mathfrak Z_{n-2}} + \alpha_1 \|\log |x|\|_\infty \|x\|.
\end{eqnarray*}
	But, also,
	$$\bigstar \geq \left\| x\left[\frac{2^j}{j!}\log^j |x|\right]_{j=2}^{n-1} - \alpha_1 \log |x| \KP_{1,n-2}(x)\right\|_{\mathfrak Z_{n-2}}  - \alpha_1 \|\log |x|\|_\infty \|x\|.$$
	
	In other words, up to a term of order $\|\log |x|\|$,
	$$\bigstar \sim \left\| x\left[\frac{2^j}{j!}\log^j |x|\right]_{j=n-1}^{2} - \alpha_1 \log |x| \KP_{1,n-2}(x)\right\|_{\mathfrak Z_{n-2}} .$$
	Now, we make a sleight of hand and from now on we will call $\bigstar$ the right term above. 	
	\begin{eqnarray*}
		\bigstar&=& \left\| x\left[\frac{2^j}{j!}\log^j |x|\right]_{j=n-1}^{2} - x \alpha_1 \log |x|\left[\frac{2^j}{j!}\log^j |x|\right]_{j=n-1}^{1}\right\|_{\mathfrak Z_{n-2}} \\
		&=& \left\| x \quad \left[\frac{2^j}{j!}\log^j |x|\right]_{j=n-1}^{2} - \alpha_1 \left[\frac{2^j}{j!}\log^{j+1} |x|\right]_{j=n-2}^{1}\right\|_{\mathfrak Z_{n-2}} \\
		&=& \left\| x \quad \left[\frac{2^{j+1}}{(j+1)!}\log^{j+1} |x|\right]_{j=n-2}^{1} - \alpha_1 \left[\frac{2^j}{j!}\log^{j+1} |x|\right]_{j=n-2}^{1}\right\|_{\mathfrak Z_{n-2}} \\
		&=& \left\| x \quad \left[\left( \frac{2^{j+1}}{(j+1)!} - \alpha_1 \frac{2^j}{j!}\right) \log^{j+1} |x|\right]_{n-2}^{1}\right\|_{\mathfrak Z_{n-2}} \end{eqnarray*}
	with rightmost coefficient $\alpha_2 = \frac{2^2}{2!} - \alpha_1 \frac{2}{1!}$.\\
	
	Step 3. Keeping up with the same ideas
	\begin{eqnarray*}
		\bigstar&=& \left\| x \quad \left[\left( \frac{2^{j+1}}{(j+1)!} - \alpha_1 \frac{2^j}{j!}\right) \log^{j+1} |x|\right]_{j=n-2}^{1}\right\|_{\mathfrak Z_{n-2}} \\
		&\sim& \left\| x \quad \left[\left( \frac{2^{j+1}}{(j+1)!} - \alpha_1 \frac{2^j}{j!}\right) \log^{j+1} |x|\right]_{j=n-2}^{2}
		- \alpha_2\log^2 |x| \left [\frac{2^j}{j!}\log^j |x| \right]_{j=n-3}^{1}\right\|_{\mathfrak Z_{n-3}} \end{eqnarray*}
	up to a term of order $\alpha_2\|\log^2 |x|\|_\infty $. Thus, after a new sleight of hand
	\begin{eqnarray*}
		\bigstar&=& \left\| x \quad \left[\left( \frac{2^{j+1}}{(j+1)!} - \alpha_1 \frac{2^j}{j!}\right) \log^{j+1} |x|\right]_{j=n-2}^{2}
		- \alpha_2 \left [\frac{2^j}{j!}\log^{j+2} |x| \right]_{j=n-3}^{1}\right\|_{\mathfrak Z_{n-3}} \\
		&=&\left\| x \quad \left[\left( \frac{2^{j+2}}{(j+2)!} - \alpha_1 \frac{2^{j+1}}{(j+1)!}\right) \log^{j+2} |x| - \alpha_2 \frac{2^j}{j!}\log^{j+2} |x| \right]_{j=1}^{n-3}\right\|_{\mathfrak Z_{n-3}} \\
		&=&\left\| x \quad \left[\left( \frac{2^{j+2}}{(j+2)!} - \alpha_1 \frac{2^{j+1}}{j+1!} - \alpha_2 \frac{2^j}{j!}\right) \log^{j+2} |x| \right]_{j=n-3}^{1}\right\|_{\mathfrak Z_{n-3}}
	\end{eqnarray*}
	The pattern is clear now, and we get that, after $n-1$ steps,
	$$\bigstar \sim \left\| x \quad \left(\frac{2^{n-1}}{(n-1)!} - \alpha_1 \frac{2^{n-2}}{(n-2)!} - \cdots - \alpha_{n-1} \frac{2}{1!} \right) \log^{n-1} |x|\right\|.$$
Therefore, to making this term finite one has $x\log^{n-1} |x|\in \ell_2$, namely $\sum x^2 \log^{2(n-1)} | x|$ is finite; and this means that $x$ belongs to the Orlicz function space $\ell_{f_{n-1}}$.\end{proof}

The following lemma is the tool to translate results from the Rochberg spaces to Orlicz spaces.

\begin{lemma}\label{L_domain} If $(u_v)_{v\in \N}$ is a normalized block basic sequence in $\ell_2$ then $w_v=(\KP_{1,n-1}(u_v),u_v)$ is equivalent to the canonical basis of $\Dom(\KP_{1,n-1})$.\end{lemma}
\begin{proof} We shall prove that if $(x_v)_{v\in\mathbb{N}}$ is a sequence of scalars in $\ell_2$ then
		$$\sum_vx_vw_v\in\mathfrak Z_n\quad\text{if and only if}\quad\sum_vx_ve_v\in\Dom\KP_{1,n-1}.$$ In \cite[Theorem 7.3]{gspaces} it is proved that \begin{equation*}
		T_U=\begin{pmatrix}
		\mathfrak u & \KP^1 \mathfrak{u}  & \KP^2 \mathfrak{u}  & \cdots & \KP^{n-1} \mathfrak{u} \\
		0 & \mathfrak u & \KP^1 \mathfrak{u} & \KP^2 \mathfrak{u}   & \cdots\\
		0 & 0 & \mathfrak u & \KP^1 \mathfrak{u}  & \KP^2 \mathfrak{u}\\
		0 & 0 & 0 & \mathfrak u & \KP^1 \mathfrak{u}\\
		0 & 0 & 0 & 0 & \mathfrak u &\\
		\end{pmatrix}
		\end{equation*}
		where $\KP^k \mathfrak{u} =\frac{2^k}{k!}u_k\log^k|u_k|$ defines a linear into isomorphism on $\mathfrak R_n$ (when $k=0$ this is, in particular, the usual multiplication operator $\mathfrak{u}(e_v)=u_v$ on $\ell_2$). Observe that
		\begin{eqnarray*}
T_U(0,\ldots,0,\sum_vx_ve_v)&=&\Big(\KP_{n-1}(\mathfrak{u})(\sum_vx_ve_v),\ldots,\KP_1(\mathfrak{u})(\sum_vx_ve_v),\mathfrak{u}(\sum_vx_ve_v) \Big)\\
		&=&\Big(\sum_vx_v\KP_{n-1}(\mathfrak{u})(e_v),\ldots,\sum_vx_n\KP_1(\mathfrak{u})(e_v),\sum_vx_v\mathfrak{u}(e_v)\Big)\\
		&=&\Big(\sum_vx_v\KP_{n-1}(u_v),\ldots,\sum_vx_v\KP_1(u_v),\sum_vx_vu_v\Big)\\
		&=&\Big(\sum_vx_v\big(\KP_{n-1}(u_v),\ldots,\KP_1(u_v)\big),\sum_vx_vu_v\Big)\\
		&=&\Big(\sum_vx_v\KP_{1,n-1}(u_v),\sum_vx_vu_v\Big)\\
		&=&\sum_vx_vw_v.\end{eqnarray*}
Thus, $\|\sum_vx_vw_v\|=\|T_U(0,\ldots,0,\sum_vx_ve_v)\|\sim \|(0,\ldots,0,\sum_vx_ve_v)\|_{\mathfrak Z_n}$ yields that $\sum_vx_vw_v$ converges if and only if $(0,\ldots,0,\sum_vx_ve_v)$ converges in $\mathfrak Z_n$; equivalently, if and only if $\sum_vx_ve_v$ converges in $\Dom\KP_{1,n-1}$.\end{proof}

We show now that the $n^{th}$ Rochberg space has exactly $n$ types of basic sequences.

\begin{theorem}\label{seqinZ3}\label{T_subsequence}
Every normalized basic sequence in $\mathfrak Z_n$ admits a subsequence equivalent to the basis of one of the spaces $\ell_{f_j}$, $1\leq j\leq n$.\end{theorem}
\begin{proof} We proceed by induction, recalling that the initial step was proved by Kalton and Peck \cite[Theorem 5.4]{kaltpeck}. Let $(x_{n-1}^k, \dots, x_0^k)_{k\in \N}$ be a normalized basic sequence in $\mathfrak Z_n$. If $\|x_0^k\|\To 0$  we can assume that $\sum \|x_0^k\|<\infty$ and thus that, up to a perturbation, $(x_{n-1}^k, \dots, x_1^k)$  is a basic sequence in $\mathfrak Z_{n-1}$. By the induction hypothesis, this sequence admits a subsequence equivalent to the basis of one of the spaces $\ell_{f_j}$, $1\leq j\leq n-1$.

If $\|x_0^n\|\geq \varepsilon$ then we can assume after perturbation that there is a block basic sequence $(u_k)_{k\in\mathbb{N}}$ in $\ell_2$ so that $\sum\|x_0^k-u_k\|<\infty$. Since
	\begin{align*}
	(x_{n-1}^k,\ldots,x_0^k)&=(x_{n-1}^k,\ldots,x_0^k)-(\KP_{1,n-1}(u_k),u_k)+(\KP_{1,n-1}(u_k),u_k)\\
	&=\big((x_{n-1}^k,\ldots,x_1^k)-\KP_{1,n-1}(u_k),x_0^k-u_k\big)+(\KP_{1,n-1}(u_k),u_k).
	\end{align*}
	and $x_0^k-u_k\to0$, we can assume (by the first part of the proof) that the sequence $\big((x_{n-1}^k,\ldots,x_1^k)-\KP_{1,n-1}(u_k),x_0^k-u_k\big)_k$ admits a subsequence equivalent to the canonical basis of $\ell_{f_j}$ for some $0\leq j\leq n-2$. Thus, Proposition \ref{basisorlicz} plus Lemma \ref{L_domain} yield that the sequence $\big(\KP_{1,n-1}(u_k),u_k\big)_k$ is equivalent to the canonical basis of $\ell_{f_{n-1}}$. We conclude that a subsequence $(w_k)_k$ of our starting $(x_{n-1}^k,...,x_0^k)_k$ is equivalent to the canonical basis of $\ell_{f_{n-1}}$: indeed, if $\sum t_k \big(\KP_{1,n-1}(u_k),u_k\big)$ converges then $\sum t_k \big((x_{n-1}^k,\ldots,x_1^k)-\KP_{1,n-1}(u_k),x_0^k-u_k\big)$ converges because $\ell_{f_{n-1}}\subset \ell_{f_{n-2}}\subset\cdots\subset\ell_{f_1}\subset\ell_2$. Passing to a subsequence, the sum
	$$\sum t_k (x_{n-1}^k,\ldots,x_0^k)=\sum_k t_k\big((x_{n-1}^k,\ldots,x_1^k)-\KP_{1,n-1}(u_k),x_0^k-u_k\big)+\sum t_k\big(\KP_{1,n-1}(u_k),u_k\big)$$
	converges. Conversely, if $\sum t_kw_k$ converges, then $\sum t_k\big(\KP_{1,n-1}(u_k),u_k\big)$ converges. If not, then
	\begin{eqnarray*}
	\Big\|\sum t_kw_k\Big\|_{\mathfrak{R}_n}&=&\Big\|\sum t_k\big((x_{n-1}^k,\ldots,x_1^k)-\KP_{1,n-1}(u_k),x_0^k-u_k\big)+\sum t_k\big(\KP_{1,n-1}(u_k),u_k\big)\Big\|_{\mathfrak Z_n}\\
	&\geq& \Big\|\sum t_k\big(\KP_{1,n-1}(u_k),u_k\big)\Big\|_{\mathfrak Z_n}- \Big\|\sum t_k\big((x_{n-1}^k,\ldots,x_1^k)-\KP_{1,n-1}(u_k),x_0^k-u_k\big)\Big\|_{\mathfrak Z_n}\\
	&\sim& \Big\|\sum t_ke_k\Big\|_{\ell_{f_{n-1}}}-\Big\|\sum t_k e_k\Big\|_{\ell_{f_j}},
	\end{eqnarray*}
and since $f_{n-1}\gg f_j$ it follows that $\sum_k t_kw_k$ does not converge, which is a contradiction. \end{proof}

Recall that a Banach space $X$ has GL.l.u.st (Local Unconditional Structure) if there exists a constant $K$ such that for every finite dimensional Banach space $F\subset X$, the inclusion map $i_F:F\hookrightarrow X$ factorizes through a Banach space $Y$ with unconditional basis so that $\|U\|\,\|V\|\,uc(Y)\leq K$, where $UV=i_F$ and $uc(Y)$ is the unconditionality constant of $Y$.

\begin{theorem}\label{lattice} For $n>1$ the space $\mathfrak Z_n$ does not contain complemented subspaces isomorphic to a Banach lattice (resp. with GL.l.u.st) and is not a complemented subspace of a Banach space isomorphic to a Banach lattice (resp. with GL.l.u.st).
\end{theorem}
\begin{proof} We prove first the assertion about GL.l.u.st: Casazza and Kalton proved in \cite[Th. 3.8]{CasKal} that a Banach space that has both GL.l.u.st and a $k$-UFDD must have an unconditional basis. The Rochberg space $\mathfrak Z_n$ has an $n$-UFDD formed by the $n$-dimensional subspaces $(E_m)$ with $E_m = {\mathrm span} \{(e_m,0\ldots,0),(0,e_m,0\ldots,0),\ldots,(0,\ldots,0,e_m)\}$ so we need to prove that $\mathfrak Z_n$ does not contain complemented subspaces with unconditional basis.

Recall from \cite[Corollary 13]{katir:98} that if $M$ is an Orlicz function satisfying the $\Delta_2$-condition then the space $\ell_M$ has cotype $2$ if and only if there exists $K>0$ such that $M(tx)\geq Kt^2M(x)$ for all $0\leq t,x\leq 1$. From this test and the previous description of $\ell_{f_j}$ for each $0\leq j\leq n-1$ it follows that all such Orlicz spaces have cotype $2$ and their dual spaces $\ell_{f_j}^*$ have type $2$. Consequently \cite[Proposition 7.3]{sym} every copy of $\ell_2$ in $\ell_{f_j}$ contains a complemented copy of $\ell_2$. Thus, if $(x_n)$ were an unconditional basic sequence in $\mathfrak Z_n$ generating a complemented subspace, it would admit a subsequence  $(x_{n_k})$ equivalent to the basis of one of the spaces $\ell_{f_j}$ for $j=0,...,n-1$ by Theorem \ref{seqinZ3}. This subsequence would generate a new complemented subspace of $\mathfrak Z_n$ and therefore we conclude that $\mathfrak Z_n$ contains a complemented copy of $\ell_2$, which is impossible by Proposition \ref{seven}. Since Banach lattices have GL.l.u.st and complemented subspaces of Banach spaces with GL.l.u.st also have GL.l.u.st, and any complemented subspace of $\mathfrak Z_n$ contains a complemented copy of $\mathfrak Z_n$ by Proposition \ref{seven} the proof of the theorem is complete.\end{proof}

A different proof, inspired by the original proof \cite{JLS_GL.l.u.st} for $Z_2$ can be followed in \cite{cp}. Observe that all the spaces $\mathfrak Z_n$ have the GL-property (i.e., every absolutely summing operator $\mathfrak Z_n\To X$ factors through some $L_1$-space) according to \cite[Proposition 1]{JLS_GL.l.u.st}. The previous proof also yields that if $\mathfrak Z_n(L_2)$ denote the $n$-th Rochberg space associated to the interpolation pair $(L_1,L_\infty)$ at $1/2$
then $\mathfrak Z_n(L_2)$ are not Banach lattices since any disjoint sequence of functions in such space spans a complemented copy of $\mathfrak Z_n$. In particular, they are not K\"othe function spaces (this remark is due to F\'elix Cabello).

\section{Operator ideals in $\mathfrak Z_n$, with applications}\label{rochberg}

 Following Pietsch \cite{pietsch}, an \emph{operator ideal} $\mathcal{A}$ is a subclass of the class $\mathfrak L$ of bounded operators between  Banach spaces such that finite range operators belong to $\mathcal{A}$, $\mathcal{A}+\mathcal{A} \subset \mathcal{A}$ and {$\mathfrak L\,\mathcal{A}\,\mathfrak L\subset\mathcal{A}$.} An operator is called \emph{strictly singular} if its restriction to any infinite dimensional subspace is not an isomorphism; it is called \emph{strictly cosingular} if its composition with any quotient map by an infinite codimensional subspace is not onto. The classes of strictly singular and strictly cosingular operators, denoted respectively by $\mathfrak S$ and $\mathfrak C$, are operator ideals \cite[1.9 and 1.10]{pietsch}. Let $T\in \mathfrak L(X,Y)$.
Then $T$ is \emph{upper semi-Fredholm,} $T\in \Phi_+(X,Y)$, if has closed range and finite dimensional kernel; $T$ is \emph{lower semi-Fredholm,} $T\in \Phi_-(X,Y)$, if its range is finite codimensional, and  $\Phi=\Phi_+\cap\Phi_-$ is the class of Fredholm operators. $T\in \mathfrak L(X,Y)$ is \emph{inessential}, denoted $T\in \mathfrak{In}$, if $I_X-AT$ is a Fredholm operator for all $A\in\mathfrak L(Y,X)$ or, equivalently, $I_Y-TA$ is Fredholm for all $A\in\mathfrak L(Y,X)$. Inessential operators were introduced by Kleinecke \cite{Kleinecke}, who proved that $\mathfrak {In}$ is a closed operator ideal containing both $\mathfrak S$ and $\mathfrak C$.
%
Standard well-known facts are:
\begin{itemize}
\item The quotient map $\pi_{n,m}: \mathfrak Z_n \To \mathfrak Z_m$ is strictly singular \cite[Proposition 9]{cck}.
\item The embedding $\imath_{n,m}: \mathfrak Z_n \To \mathfrak Z_m$ is strictly cosingular (plain duality).
\end{itemize}

We need to recall here the transformation  $T\to T^+$ from \cite{symplectic} and the notion of block operator from
\cite{symplectic,gspaces}.

\adef Given an operator $T\in \mathfrak L(\mathfrak Z_n)$ we define the operator $T^+\in \mathfrak L(\mathfrak Z_n)$ by the identity $\langle T^+x,y\rangle =\langle x, Ty \rangle$.\zdef

By Proposition \ref{dualRochberg} we deduce that there exist a commutative diagram relating the dual operator $T^*$ with $T^+$ in the following sense:
\begin{equation}\label{DiagramaDual}
\xymatrixcolsep{1.5cm}\xymatrixrowsep{1.5cm}\xymatrix{\mathfrak Z_n \ar[r]^{D_n} \ar[d]^{T^+} & \mathfrak Z_n^* \ar[d]^{T^*} \\
	\mathfrak Z_n \ar[r]^{D_n} &  \mathfrak Z_n^*
}
\end{equation}

\adef Let $u$ be a disjointly supported sequence of normalized blocks in $\ell_2$. We inductively define the block operator $T_U^n\in \mathfrak L(\mathfrak Z_n)$ as: $T_U^1=u$, and for $n>1$, $T_U^n$ is the upper triangular operator making the diagrams

$$\xymatrix{0\ar[r]&\mathfrak Z_n \ar[d]_{T_U^n} \ar[r]&\mathfrak Z_{n+m} \ar[d]^{T_U^{n+m}} \ar[r]& \mathfrak Z_{m}\ar[r]\ar[d]^{T_U^m}&0\\
0\ar[r]&\mathfrak Z_n \ar[r]&\mathfrak Z_{n+m} \ar[r]& \mathfrak Z_{m}\ar[r]&0}$$ commutative
\zdef
In \cite[Proposition 7.1]{gspaces} the reader can find an explicit description of these operators. The key property of block operators proved in \cite[Prop. 6.2]{symplectic} is that $\omega_n(T_U^nx,T_U^ny)=\omega_n(x,y)$ for any $x,y\in\mathfrak{Z}_n$, or equivalently, that $(T_U^n)^+T_U^n=I_{\mathfrak{Z_n}}$. Thus is not hard to show that block operators are into isometries with complemented range (see \cite{symplectic}). Our purpose now is to show that operators on $\mathfrak Z_n$ behave in an analogous way to operators in $\mathfrak Z_2= Z_2$, as shown in \cite{kaltsym,cgp}. We begin extending the following technical result.

\begin{prop}[Essentially Kalton \cite{kaltsym}]\label{C_01}
Let $T:\mathfrak Z_n\rightarrow\mathfrak Z_n$ be any operator. If $T^+T$ is strictly singular then $T$ is strictly singular.\end{prop}
\begin{proof} We need the following result from \cite{symplectic}:

\begin{lemma}\label{P_02} If $T:\mathfrak Z_n\rightarrow\mathfrak Z_n$ is not strictly singular then there exists $\alpha\neq0$ and block operators $T_U^n$ and $T_V^n$ such that $TT_W^n-\alpha T_V^n$ is strictly singular.\end{lemma}

The proof is now clear: if $T$ is not strictly singular there exists $\alpha\neq0$ and block operators $T_U^n,T_V^n$ such that $ TT_U^n=\alpha T_V^n-S$ with $S$ strictly singular. Therefore
	$$(T_U^n)^+T^+TT_U^n=(TT_U^n)^+TT_U=(\alpha (T_V^n)^+-S^+)(\alpha T_V-S)=\alpha'(T_V^n)^+T_V^n+S'=\alpha'I_{\mathfrak Z_n}+S'.$$
	Since $T^+T$ is strictly singular, this implies that $I$ is strictly singular, a contradiction.\end{proof}

We extend now two important results of Kalton about the behaviour of operators in $Z_2$ \cite{kaltsym} to higher order Rochberg spaces (see the first part for $n=3$ in \cite{sym}). We have:

\begin{prop}\label{seven}$\;$
\begin{enumerate}
\item \label{T_01} Every operator $\tau:\mathfrak Z_n\rightarrow X$ is either strictly singular or an isomorphism on a complemented copy of $\mathfrak Z_n$.
\item \label{T_compl} Every operator $\tau: \mathfrak Z_n\rightarrow\mathfrak Z_n$ is either strictly singular or an isomorphism on a complemented subspace $E\cong\mathfrak Z_n$ such that $\tau [E]$ is also complemented.\end{enumerate}
\end{prop}

\begin{proof} We prove (1). Consider the representation $\xymatrix{0\ar[r]&\ell_2\ar[r]^{\imath_{1,n}}&\mathfrak Z_n\ar[r]^{\pi_{n,1}}&\mathfrak Z_{n-1}\ar[r]&0}$ and assume that $\tau$ is not strictly singular.  We consider two cases:
\begin{itemize}
\item Case 1. Suppose that $\tau_{| \ell_2}$ is an embedding. Then the result is obtained with the same proof of \cite[Proposition 7.6]{sym} by changing $\mathfrak Z_3$ with $\mathfrak Z_n$.
\item Case 2. Suppose that $\tau_{|Y}$ is an isomorphism for an infinite dimensional closed subspace $Y$ of $\mathfrak Z_n$. Since $\pi_{n, n-1}$ is strictly singular, there exists a basic sequence $(v_n) \subset  Y$ such that for some block basic sequence $(u_n)\subset  \ell_2$ we have $\lim \|\imath_{1,n}(u_n)- v_n\|=0$. It follows that  ${\tau T^n_{U}}_{|\ell_2}$ is an embedding for some subsequence $U$ of $(u_n)_n$.  Then the results follows from Case 1 applied to $\tau T^n_{U}$ and the fact that the range of $T^n_U$ is complemented.\end{itemize}

We prove (2). Suppose that $\tau$ is not strictly singular. By Lemma \ref{P_02} there exist block operators $T_U^n$, $T_V^n$ and $\alpha\neq0$ such that $\tau T_U^n=\alpha T_V^n+S$ with $S$ strictly singular. Thus, it follows that
$\alpha^{-1}(T_V^n)^+\tau T_U^n=\alpha^{-1}\alpha (T_V^n)^+T_V^n+(T_V^n)^+S=I+S_1.$ Thus $\alpha^{-1}(T_V^n)^+\tau T_U^n$ is a Fredholm operator with index $0$ by \cite[2.c.10]{lindtzaf1}. Let $F$ be a (closed) complement of the kernel of $(T_V^n)^+\tau T_U^n$, and let $E=T_U^n(F)$. Then $\tau$ is an isomorphism on $E$, the subspaces $E$ and $\tau(E)$ are both isomorphic to $\mathfrak Z_n$, and both of them are complemented in $\mathfrak Z_n$. Indeed, we can write $\mathfrak Z_n=(T_V^n)^+\tau T_U^n(F)\oplus N$ with $\dim N<\infty$; hence $\mathfrak Z_n=\tau(E)\oplus \big((T_V^n)^+\big)^{-1}(N)=E\oplus \big((T_V^n)^+\tau\big)^{-1}(N)$.\end{proof}

\begin{theorem}\label{singcosingr} $\mathfrak S(\mathfrak Z_n) = \mathfrak C(\mathfrak Z_n)$ is the only nontrivial maximal ideal of $\mathfrak L(\mathfrak Z_n)$.\end{theorem}

\begin{proof} Let us show first that strictly singular and cosingular operators on $\mathfrak Z_n$ coincide. For this end suppose that $T\in \mathfrak{S}(\mathfrak Z_n)$. Then $T^*\in\mathfrak{C}(\mathfrak Z_n)$ and by (\ref{DiagramaDual}) we deduce that $T^+\in\mathfrak{C}(\mathfrak Z_n)$. Moreover, by Theorem \ref{seven}, $T^+$ is strictly singular or an isomorphism on a complemented copy $E$ of $\mathfrak Z_n$ with $T^+[E]$ complemented. In the second case we would easily conclude $Id\in\mathfrak{C}(\mathfrak Z_n)$, which is false. Then $T^*$ is strictly singular, hence $T$ is strictly cosingular. If $T:\mathfrak Z_n\rightarrow\mathfrak Z_n$ is strictly cosingular then $T^*:\mathfrak Z_n^*\rightarrow\mathfrak Z_n^*$ is strictly singular. Thus $T^+$ is strictly singular, which implies $T^+T$ strictly singular, hence so is $T$ by Proposition \ref{C_01}. We show now that every proper ideal in $\mathfrak{L}(\mathfrak Z_n)$ is contained in $\mathfrak{S}(\mathfrak Z_n)$. Suppose that $T\in\mathfrak{L}(\mathfrak Z_n)$ is not strictly singular. By Theorem \ref{seven}, $T$ is an isomorphism on subspace $E$ isomorphic to $\mathfrak Z_n$ so that both $E$ and $T[E]$ are complemented. Thus there exist operators $A,B\in\mathfrak{L}(\mathfrak Z_n)$ such that $ATB=I_{\mathfrak Z_n}$. Hence $T$ does not belong to any proper ideal in $\mathfrak{L}(\mathfrak Z_n)$.\end{proof}

Diagram (\ref{DiagramaDual}) says that $T^+=D_n^{-1}T^* D_n\in \mathfrak{L}(\mathfrak Z_n)$, and since $D_n$ is an isomorphism, most of the properties of $T^+$ coincide with those of $T^*$. Namely, $\|T\|\sim \|T^+\|$, $R(T)$ is closed if and only if $R(T^+)$ is so, $T$ is an isomorphism into if and only if $T^+$ is surjective, and $T\in\Phi_-$ if and only if $T^+\in\Phi_+$. In particular, $T^+=T$ and $T\in\Phi_+$ implies $T\in\Phi$. We now show that the strictly singular and the upper semi-Fredholm character of $T:\mathfrak Z_n \rightarrow\mathfrak Z_n$ only depends on its behaviour on the natural copy $\ell_2\subset\mathfrak{R_n}$.

\begin{lemma}\label{restrict}
Let $T\in \mathfrak{L}(\mathfrak Z_n)$.
\begin{itemize}
\item[(a)] If $T\imath_{1,n}\in\Phi_+$ then $T\in\Phi_+$.
\item[(b)] If $T\in\Phi_+$  then $T^+T\in\Phi$.
\end{itemize}
\end{lemma}
\begin{proof} (a) Suppose that $T\not\in\Phi_+$. Then there exists an infinite dimensional subspace $M \subset \mathfrak Z_n$ such that $T|_M$ is an isomorphism. Since $\pi_{n,n-1}$ is strictly singular, a perturbation argument provides an infinite dimensional subspace $N\subset M$ and a compact operator $K\in \mathfrak{L}(\mathfrak Z_n)$ with $\|K\|<1$ such that $(I-K)N\subset \imath_{1,n}[\ell_2]= N(\pi_{n,n-1})$.
We shall show that $T\imath_{1,n}|_{I-K)[N]}$ is compact, hence $T\imath_{1,n} \not\in\Phi_+$. Indeed, if $(x_n)$ is a bounded sequence in $N$, then $(T(I-K)x_n)= (Tx_n-TKx_n)$ has a convergent subsequence. Since $TK$ is compact, $(Tx_n)$ has a convergent subsequence. We prove (c). Since $(T^+T)^+=T^+T$, it is enough to show that $T^+T\in\Phi_+$. Suppose that $T^+T\not\in\Phi_+$. Then $T^+T\imath_{1,n}\not\in\Phi_+$; hence there exists a normalized block basis sequence $(w_n)$ in $\ell_2$ such that, if we denote by let $i_w: \ell_2\to \ell_2$ the isometric embedding defined by $i_w e_n= w_n$, then $T^+T\imath_{1,n} i_w$ is compact. Let $T_W^n$ be the block operator associated to the sequence $(w_n)$. Since $T^+TT_W^n\imath_{1,n}e_n= T^+Ti_{1,n} i_we_n$ for each $n\in\N$, $T^+TT_W^ni_{1,n}$ is compact; hence $T^+TT_W^n$ is strictly singular by Proposition \ref{restrict}.
Thus $(T_W^n)^+T^+TT_W^n\in\mathfrak S$, hence $TT_W^n\in \mathfrak S$, so $T\notin\Phi_+$.\end{proof}

\begin{lemma}\label{L_type_sequence} Let $(x_n)_n$ be a seminormalized basic sequence sequence equivalent to the canonical basis of $\ell_{f_k}$ and let $S\in\mathfrak S(\mathfrak Z_n)$. Then $(Sx_n)_n$ has either a norm null subsequence or a seminormalized basic subsequence equivalent to the canonical basis of $\ell_{f_j}$ for some $j<k$.\end{lemma}
\begin{proof} Assume that $(Sx_n)_n$ has no norm-null subsequence. Then passing to a subsequence que have that $\|Sx_n\|\geq\varepsilon>0$. Thus, $(Sx_n)_n$ is seminormalized and weakly null since $S$ is bounded and $(x_n)_n$ is basic in a reflexive space, hence shrinking. Using Bessaga-Pe\l czy\'{n}ski selection principle we can assume that $(Sx_n)_n$ is seminormalized and basic. By Theorem \ref{T_subsequence} we have that $(Sx_n)_n$ has a (seminormalized basic) subsequence $(Sx_{n_l})_l$ equivalent to the canonical basis of $\ell_{f_j}$ for some $0\leq j\leq k-1$. Indeed, $(Sx_{n_l})_l$ can not be equivalent to the canonical basis of $\ell_{f_k}$ because $S$ is strictly singular. Moreover, we have the chain of continuous and proper inclusions
		\begin{equation}\label{eq_chain_Orlicz}
		\ell_{f_{n-1}}\subset\ell_{f_{n-2}}\subset\cdots\subset \ell_{f_1}\subset\ell_2
		\end{equation}
		and those spaces are not isomorphic. Then, if $(Sx_{n_l})_l$ were equivalent to the canonical basis of $\ell_{f_j}$ for some $j>k$ then, up to a subsequence,
		$$\|\sum a_l e_l\|_{\ell_{f_{j}}}\sim \|\sum a_lSx_{n_l}\|\leq \|S\|\|\sum a_l x_{n_l}\|\sim \|\sum a_l e_l\|_{\ell_{f_{k}}}.$$
		Since by (\ref{eq_chain_Orlicz}) we have the inequality $\|x\|_{\ell_{f_{k}}}\leq \|x\|_{\ell_{f_j}}$, we conclude that the bases of $\ell_{f_{k}}$ and $\ell_{f_j}$ would bee equivalent, which is false.\end{proof}

In $\ell_2$ compact and strictly singular operators coincide. In \cite{cgp} we proved that the composition of two strictly singular operators is compact. We generalize this fact:

\begin{theorem}\label{composition}	If $S_1,\dots, S_n\in \mathfrak S(\mathfrak Z_n)$ then $S_1\cdots S_n \in \mathfrak K(\mathfrak Z_n)$.\end{theorem}
\begin{proof} Recall that an operator $T$ on a reflexive Banach space is compact if for every normalized weakly null sequence $(x_n)_n$, we have that $(Tx_n)_n$ has a norm-null subsequence. Thus, let $(x_n)_n$ be a normalized weakly null sequence. Then passing to a subsequence and using Bessaga-Pe\l czy\'{n}ski selection principle we can assume by Theorem \ref{seqinZ3} that $(x_n)_n$ is equivalent to the canonical basis of one of the spaces $\ell_{f_k}$ for some $0\leq k\leq n-1$. Consider now the sequence $(S_1x_n)_n$. By Lemma \ref{L_type_sequence} we have two possibilities:
		\begin{itemize}
			\item[(1)]{$\|S_1x_{n_k}\|\to 0$ for some subsequence;}
			\item[(2)]{$(S_1x_n)_n$ has a seminormalized basic subsequence equivalent to the canonical basis of the space $\ell_{f_j}$ for some $j<k$.}
		\end{itemize}
		In the first case $S_1$ is compact an thus the product $S_1\cdots S_n$ is also compact. In the second case we can consider the sequence $(S_2S_1x_n)_n$ and use Lemma \ref{L_type_sequence} again to reach either a norm-null subsequence or a further subsequence of $(S_2S_1x_n)_n$ that is equivalent to the canonical basis of $\ell_{f_m}$ for some $m<j<k$. It is clear that this process can be repeated at most $n$ times, the number of possible non-equivalent basic sequences, until reaching to a norm null subsequence. Thus $(S_1\cdots S_nx_n)_n$ has a norm-null subsequence, hence $S_1\cdots S_n\in\mathfrak K(\mathfrak Z_n)$. \end{proof}

We conclude observing if $T = \imath_{n-1,n}\pi_{n,n-1}: \mathfrak Z_n \To \mathfrak Z_n$ and $1\leq k\leq n$ then $T^{k}=\imath_{n-k,n}\pi_{n,n-k}$; hence, $T^k$ is a strictly singular non-compact operator for $k\leq n-1$.

\section{Complemented subspaces of $\mathfrak Z_n$}

From Proposition \ref{seven} we deduce that every infinite dimensional complemented subspace of $\mathfrak Z_n$ contains a further complemented subspace isomorphic to $\mathfrak Z_n$. The celebrated still open hyperplane problem is determining whether $Z_2$ is isomorphic to its hyperplanes. An equivalent formulation is: do there exist Fredholm operators in $\mathfrak L(Z_2)$ with odd Fredholm index? Of course that
it is not known whether $\mathfrak Z_n$ is isomorphic to its hyperplanes ($n>1$). We conjecture that Fredholm operators in $\mathfrak L(\mathfrak Z_n)$ have Fredholm index multiple of $n$ or, equivalently, $\mathfrak Z_n$ is not isomorphic to its subspaces of codimension $1,2, \dots, n-1$.
Let us add now a new piece of knowledge:

\begin{theorem}\label{T_complementado} Every subspace of $\mathfrak Z_n$ isomorphic to $\mathfrak Z_n$ is complemented.
\end{theorem}
\begin{proof}
Let $T\in \mathfrak{L}(\mathfrak Z_n)$ an isomorphism into with $R(T)=M$. Then $T^+T\in\Phi$, hence here exists a finite codimensional subspace $N$ of $R(T)$ such that $T^+|_N$ is an isomorphism and $T(N)$ is finite codimensional, hence complemented; thus $N$ is complemented and so is $M$.\end{proof}

\begin{cor}
Every semi-Fredholm operator on $\mathfrak Z_n$ has complemented kernel and range.\end{cor}
\begin{proof}
If $T\in \Phi_+(\mathfrak Z_n)$ then the kernel is finite dimensional, and the range $R(T)$ is complemented using the arguments of Theorem \ref{T_complementado}. If $T\in\Phi_-(\mathfrak Z_n)$ then $R(T)$ is finite codimensional and $T^*\in\Phi_+(\mathfrak Z_n^*)$. Taking into account that $\mathfrak Z_n\equiv\mathfrak Z_n^*$ we conclude that $R(T^*)$ is complemented using the first part of the theorem, hence $N(T)={}^\perp R(T)$ is also complemented.\end{proof}

Recall from \cite{hmbst} that a Banach space $X$ is said to be $Y$-automorphic if every isomorphism between
two infinite codimensional subspaces of $X$ isomorphic to $Y$ can be extended to an automorphism of $X$.
It is clear that $\ell_2$ is $\ell_2$-automorphic and that $Z_2$ is not $\ell_2$-automorphic: indeed, since $Z_2\simeq Z_2 \oplus Z_2$, an
isomorphism between the subspaces $\ell_2 \oplus 0$ and $\ell_2 \oplus \ell_2$ cannot be extended to an automorphism of $Z_2$. Surprisingly, one has:

\begin{prop}\label{automorphic} $\mathfrak Z_n$ is $\mathfrak Z_n$-automorphic.\end{prop}
\begin{proof} Since Theorem \ref{T_complementado} establishes that every copy of $\mathfrak Z_n$ in $\mathfrak Z_n$ is complemented we just need to worry about complemented copies with infinite dimensional complement. Hence, our task is to show that those complements are isomorphic. Suppose that $\mathfrak Z_n \simeq \mathfrak Z_n \oplus F$. By (2) of Proposition \ref{seven} one has $F\simeq \mathfrak Z_n \oplus N$, and therefore $F \oplus F \simeq F \oplus \mathfrak Z_n \oplus N \simeq \mathfrak Z_n \oplus N \simeq F$. It follows from Pe\l czy\'nski's decomposition argument that if $E$ is a complemented subspace of $\mathfrak Z_n$ and $E \oplus E \simeq E$ then $E$ is isomorphic to $\mathfrak Z_n$.\end{proof}

\section{The compactness theorem revisited}

In \cite[Corollary 5.9]{ccfm} it was proved a remarkable result: if one has a commutative diagram
$$\xymatrix{
0 \ar[r] & \ell_2 \ar[r]\ar[d]_{T} & Z_2 \ar[r] \ar[d] & \ell_2  \ar[r] \ar[d]^{S} & 0\\
0 \ar[r] & \ell_2  \ar[r] & Z_2 \ar[r] & \ell_2  \ar[r]  & 0}$$
then $T-S$ is compact. It is not known if something similar occurs for any other twisted Hilbert space, but we will show now that higher order Rochberg spaces admit similar results. Observe that the hypothesis ``one has a commutative diagram" can be reformulated as
``there is an upper triangular operator $Z_2\To Z_2$ of the form $\left(
                                                  \begin{array}{cc}
                                                    T & \square \\
                                                    0 & S \\
                                                  \end{array}
                                                \right)$". In general, if
$$R= \left( \begin{array}{cccccc}
a_{11}  & a_{12} & a_{13} & \cdots & a_{1,n-1} & a_{1n}\\
0 & a_{22} & a_{23} & a_{24} & \cdots & a_{2n} \\
0 & 0 & a_{33} & a_{34} & a_{35} & \cdots\\
0 & 0 & 0 & \cdots  & \cdots & \cdots\\
0 & 0 &  0 & \cdots  & a_{n-1 n-1} & a_{n-1 n}\\
0 & 0 & \cdots & 0 & 0 & a_{nn}
\end{array}\right)\in \mathfrak L(\mathfrak Z_n)$$ and $k<n$ then
$$R_k = \left(\begin{array}{cccccc}
a_{11}  & a_{12} & a_{13} & \cdots & a_{1,k-1} & a_{1k}\\
0 & a_{22} & a_{23} & a_{24} & \cdots & a_{2k} \\
0 & 0 & a_{33} & a_{34} & a_{35} & \cdots\\
0 & 0 & 0 & \cdots  & \cdots & \cdots\\
0 & 0 &  0 & \cdots  & a_{k-1 k-1} & a_{k-1 k}\\
0 & 0 & \cdots & 0 & 0 & a_{kk}
\end{array}\right) \in \mathfrak L(\mathfrak Z_k)$$
and the corresponding induced operator on the quotient space is
$$R^{n-k}= \left(\begin{array}{cccccc}
a_{k+1k+1}  & a_{k+1k+2} & a_{k+1k+3} & \cdots & a_{k+1,n-1} & a_{k+1n}\\
0 & a_{k+2k+2} & a_{k+2k+3} & a_{k+2k+4} & \cdots & a_{k+2n} \\
0 & 0 & a_{k+3k+3} & a_{k+3k+4} & a_{k+3k+5} & \cdots\\
0 & 0 & 0 & \cdots  & \cdots & \cdots\\
0 & 0 &  0 & \cdots  & a_{n-1 n-1} & a_{n-1 n}\\
0 & 0 & \cdots & 0 & 0 & a_{nn}
\end{array}\right)\in \mathfrak L(\mathfrak Z_{n-k})$$
so that there is a commutative diagram
$$\xymatrix{
0 \ar[r] & \mathfrak Z_k \ar[r]\ar[d]_{R_k} & \mathfrak Z_{n} \ar[r] \ar[d]^R & \mathfrak Z_{n-k} \ar[r] \ar@/_2pc/[ll]_{\KP_{n-k,k}}
\ar[d]^{R^{n-k}} & 0\\
0 \ar[r] & \mathfrak Z_k  \ar[r] & \mathfrak Z_{n} \ar[r] & \mathfrak Z_{n-k} \ar[r] \ar@/^2pc/[ll]^{\KP_{
n-k,k}} & 0}.$$

We obtain now the desired generalization:

\begin{theorem} If $R\in \mathfrak L(\mathfrak Z_n)$ then $R_k - R^k \in \mathfrak S(\mathfrak Z_k)$ for every $k<n$.\end{theorem}
\begin{proof} We work inductively. For $n=2$ is the result just mentioned. The case $n=3$ will help us to explain the strategy: an upper triangular operator $R = \left(
  \begin{array}{ccc}
    \alpha & \beta & \epsilon \\
    0 & \gamma & \delta \\
    0 & 0 & \eta \\
  \end{array}
\right)\in \mathfrak L(\mathfrak Z_3)$ generates two commutative diagrams:
$$\xymatrix{
0 \ar[r] & Z_2 \ar[r]\ar[d]_{R_2} & \mathfrak Z_{3} \ar[r] \ar[d]^R & \ell_2 \ar[d]^{R^1}\ar[r] & 0\\
0 \ar[r] & Z_2 \ar[r] & \mathfrak Z_{3} \ar[r] & \ell_2 \ar[r] & 0}\quad \quad \xymatrix{
0 \ar[r] & \ell_2 \ar[r]\ar[d]_{R_1} & \mathfrak Z_{3} \ar[r] \ar[d]^R & Z_2 \ar[d]^{R^2}\ar[r] & 0\\
0 \ar[r] & \ell_2 \ar[r] & \mathfrak Z_{3} \ar[r] & Z_2  \ar[r] & 0}$$
Since $R_2\in \mathfrak L(\mathfrak Z_2)$, $\alpha - \gamma$ is compact; and since $R^2\in \mathfrak L(\mathfrak Z_2)$, $(R_2)_1 - (R_2)^1 = \gamma - \eta $ is compact as well. Therefore $\alpha - \gamma$ is compact too. Since $\alpha -\gamma$ is the restriction of the operator
$$R_2 - R^2 = \left( \begin{array}{ccc}
    \alpha & \beta \\
    0 & \gamma  \\
  \end{array} \right)- \left(
  \begin{array}{ccc}
    \gamma & \delta\\
   0 & \eta \\
  \end{array}\right) = \left(
                         \begin{array}{cc}
                           \alpha - \gamma & \beta - \delta \\
                           0 & \gamma - \eta \\
                         \end{array}
                       \right)
  $$ to $\ell_2$, it follows that $R_2 - R^2$ must be strictly singular. Assume the result has been proved for $n-1$ and pick $R\in \mathfrak L(\mathfrak Z_{n})$. Using induction one obtains that $a_{ii} - a_{jj}$ is compact, and from here $R_2 - R^2$ is strictly singular too
  (as well as all ``intermediate" operators $(R_k)_2 - (R_k)^2$ and $(R^k)_2 - (R^k)^2$). From here we get that $R_3 - R^3$ is strictly singular too (as well as all ``intermediate" operators $(R_k)_3 - (R_k)^3$ and $(R^k)_3 - (R^k)^3$). And so on. The only remaining case is

  $$R_{n-1} - R^{n-1} = \left( \begin{array}{ccccc}
a_{11}  & a_{12} & a_{13} & \cdots & a_{1,n-1} \\
0 & a_{22} & a_{23} & \cdots  & a_{2,n-1} \\
0 & 0 & a_{33} & \cdots & a_{3,n-1}  \\
0 & 0 & 0 & \cdots  & \cdots  \\
0 & 0 &  0 & \cdots  & a_{n-1 n-1}\\
\end{array} \right)- \left( \begin{array}{ccccc}
a_{22} & a_{23} & a_{24} & \cdots & a_{2n} \\
0 & a_{33} & a_{34} & a_{35} & \cdots\\
0 & 0 & \cdots  & \cdots & \cdots\\
0 &  0 & \cdots  & a_{n-1 n-1} & a_{n-1 n}\\
0 & \cdots & 0 & 0 & a_{nn}
\end{array}\right)$$
but it is plain that is strictly singular since its restriction to $\ell_2$ is compact. \end{proof}

As an immediate consequence from this result we obtain:
\begin{prop}\label{compact} Given a commutative diagram
$$\xymatrix{
0 \ar[r] & \mathfrak Z_n \ar[r]\ar[d]_{T} & \mathfrak Z_{2n} \ar[r] \ar[d] & \mathfrak Z_n \ar[r]
\ar[d]^{S} & 0\\
0 \ar[r] & \mathfrak Z_n  \ar[r] & \mathfrak Z_{2n} \ar[r] & \mathfrak Z_n \ar[r] & 0}$$
in which $T$ and $S$ are upper triangular operators then $T-S$ is strictly singular.\end{prop}

We conjecture that given an upper triangular operator $T: \mathfrak Z_n \To \mathfrak Z_n$ the commutator $[T, \KP_{n,n}]$ is not singular. More precisely, that if $[T, \KP_{n,n}]$ is singular then $T$ is strictly singular.

\section{Appendix. The generalized Commutator Theorem}

In this section we  present a general version of the classical Commutator theorem that relates operators ``acting on the scale" with operators on Rochberg spaces. Commutator theorems have a long history since the classical commutator theorem, case $n=2$, was first obtained by Rochberg and Weiss \cite{rochweiss} to the general case obtained by Rochberg \cite{rochberg}; other commutator theorems have been studied in several papers, such as \cite{caceso,caceso2,racsam}. Depending the approach one tackles the commutator theorem is obvious, difficult or surprising. Its picturesque formulation is that the diagrams
$$\xymatrix{
0 \ar[r] & \mathfrak R_k \ar[r]\ar[d]_{T_k} & \mathfrak R_n \ar[r] \ar[d]^{T_n} & \mathfrak R_v \ar[r] \ar@/_2pc/[ll]_{\Omega_{v,k}}
\ar[d]^{T_v} & 0\\
0 \ar[r] & \mathfrak R_k  \ar[r] & \mathfrak R_n \ar[r] & \mathfrak R_v \ar[r] \ar@/^2pc/[ll]^{\Omega_{
v,k}} & 0}$$
are commutative. The version in \cite{racsam} is somehow the most abstract and, therefore, universally valid. The version we present next is subtly different from previous ones in that our setting introduces the (generalized) idea of operator acting on the scale in the context of Kalton spaces and deals with $\Omega_{n,k}$ for all values of $n,k$. Our proof differs from that of Rochberg and we hope it clarifies  several obscurish points in \cite{caceso2}.

A bounded operator $\tau:\Sigma \To \Sigma$ is said to act on the scale $(X_z)_z$ generated by the pair $(X_0, X_1)$ and the Calder\'on space $\mathcal C = \mathcal C(X_0, X_1)$ if $\tau[X_0]\subset X_0$ and $\tau[X_1]\subset X_1$. It is a simple observation that if $\tau$ is an operator acting on that scale then there exist an operator $T:\mathcal{C} \rightarrow\mathcal{C} $ such that for every $z$  the diagram
$$\xymatrixrowsep{1.3cm}\xymatrixcolsep{1.3cm}
\xymatrix{\mathcal{C} \ar[r]^{\Delta_0} \ar[d]_{T} & \Sigma \ar[d]^{\tau} \\
	\mathcal{C} \ar[r]^{\Delta_0} &  \Sigma}
$$commutes. The operator $T$ is necessarily defined by $T(f)(z)=\tau(f(z))$. Observe then (see \cite{racsam,caceso,sym}) that
$(Tf)'(z) = \tau (f'(z))$. Thus, if $T_n$ denotes the operator $T_n(x_1,\ldots,x_n)=(\tau x_1,\ldots, \tau x_n)$ then for any $v\geq0$ and $k>0$ such that $k+v=n$ the diagram
\begin{equation}\label{eq_interpolator}
\xymatrixrowsep{1.3cm}\xymatrixcolsep{2cm}\xymatrix{	\mathcal{C} \ar[r]^{\big(\Delta_{k+v-1},\ldots,\Delta_v\big)} \ar[d]_{T} & \Sigma^k \ar[d]^{T_k} \\
	\mathcal{C} \ar[r]^{\big(\Delta_{k+v-1},\ldots,\Delta_v\big)} &  \Sigma^k}
\end{equation}commutes. Transplanting the idea to the admissible setting can be done as follows: Let us say that a bounded operator $\tau:\Sigma \To \Sigma$ acts on the scale $(X_z)_z$ generated by the admissible space $\mathscr F$ if there exist an operator $T:\mathscr F  \rightarrow \mathscr F $ making the diagram
$$\xymatrixrowsep{1.3cm}\xymatrixcolsep{1.3cm}
\xymatrix{\mathscr F \ar[r]^{\Delta_0} \ar[d]_{T} & \Sigma \ar[d]^{\tau} \\
	\mathscr F \ar[r]^{\Delta_0} &  \Sigma}
$$commute. Thus, also if
$T_n(x_1,\ldots,x_n)=(\tau x_1,\ldots, \tau x_n)$ then for any $v\geq0$ and $k>0$ such that $k+v=n$ the diagram obtained from (\ref{eq_interpolator}) replacing $\mathcal C$ by $\mathscr F$ commutes. One has:

\begin{theorem}[Commutator Theorem]\label{Commutator}
Let $\tau$ be an operator acting on the scale $(X_z)_z$ generated by the admissible space $\mathscr F$ and $n\in\mathbb{N}$. The following equivalent statements hold:
\begin{itemize}
\item $T_n: \mathfrak R_n \To \mathfrak R_n$ is bounded.
\item For any $k,v\in\mathbb{N}$ such that $k+v=n$ $\Big\|T_k\Omega_{v,k}-\Omega_{v,k} T_v\Big\| <\infty$.
\end{itemize}\end{theorem}

\begin{proof} Let us show first that the commutator $T_k\Omega_{v,k}-\Omega_{v,k}T_v$ is bounded from $\mathfrak R_v$ into $\mathfrak R_k$. To this end, vet $B_\theta:\mathfrak R_v\rightarrow\mathscr F$ be a bounded and homogeneous selection for the quotient map $\big(\Delta_{v-1},\ldots,\Delta_1,\Delta_0\big):\mathscr F\rightarrow\mathfrak R_v$ so that $\Omega_{v,k}=\big(\Delta_{k+v-1},\ldots,\Delta_v \big)B_\theta$ and
\begin{eqnarray*}T_k\Omega_{v,k} -\Omega_{v,k}(T_v) &=& T_k\big(\Delta_{k+v-1},\ldots,\Delta_v\big)B_\theta -\big(\Delta_{k+v-1},\ldots, \Delta^v\big)B_\theta(T_v )\\
&=& \big(\Delta_{k+v-1},\ldots,\Delta_v\big)T B_\theta -\big(\Delta_{k+v-1},\ldots, \Delta^v \big)B_\theta(T_v)\\
&=&\big(\Delta_{k+v-1},\ldots,\Delta^v\big)\Big[TB_\theta -B_\theta(T_v)\Big]\end{eqnarray*}
because of the commutativity of diagram (\ref{eq_interpolator}). Again by this same commutativity,
$TB_\theta(x)-B_\theta(T_vx)\in\ker \big(\Delta_{v-1},\ldots,\Delta_1\big)$: indeed, just consider inductively the case $n=v=v'+k'$, $v'=0$ and $k'=v$. Since $\big(\Delta_{v-1},\ldots,\Delta_0\big)T=T_v\big(\Delta_{v-1},\ldots,\Delta_0\big)$ and $B_\theta$ is a selector for $\big(\Delta_{v-1},\ldots,\Delta_0\big)$, it follows that
$\big(\Delta_{v-1},\ldots,\Delta_0\big)TB_\theta -\big(\Delta_{v-1},\ldots,\Delta_0\big)B_\theta T_v =T_v\big(\Delta_{v-1},\ldots,\Delta_0\big)B_\theta -T_v =0$. Now, it has been proved in \cite[Prop. 1]{cck} that the operator $(\Delta_{k+v-1},\ldots,\Delta_v): \ker \big(\Delta_{v-1},\ldots,\Delta_1\big) \rightarrow \mathfrak R_k$ is bounded and onto for any $0<\theta<1$. Therefore $\big(\Delta_{k+v-1},\ldots,\Delta_v\big)\Big(TB_\theta -B_\theta T_v \Big)\left[ \mathfrak R_v\right] \subset \mathfrak R_k$. We conclude that the map $T_k\Omega_{v,k}-\Omega_{v,k}T_v:\mathfrak R_v\rightarrow\mathfrak R_k$ is well defined. Its boundedness follows by induction and the same type of arguments:
\begin{eqnarray*}
\Big\|T_k\Omega_{v,k}(x)-\Omega_{v,k} T_vx \Big\|_{\mathfrak R_k}&=&\Big\|T_k\big(\Delta_{k+v-1},\ldots,\Delta_v \big)B_\theta(x)-\big(\Delta_{k+v-1},\ldots,\Delta_v \big)B_\theta T_vx \Big\|_{\mathfrak R_k}\\
&=&\Big\|\big(\Delta_{k+v-1},\ldots,\Delta_v \big)TB_\theta(x)-\big(\Delta_{k+v-1},\ldots,\Delta_v\big)B_\theta T_vx\Big\|_{\mathfrak R_k}\\
&=&\|\big(\Delta_{k+v-1},\ldots,\Delta_v \big)\|_{\mathscr F}\,\Big\|TB_\theta(x)-B_\theta T_vx \Big\|_{\mathfrak R_k}\\
&\leq& C \Big(\|T\|_{\mathscr F}\|B_\theta\|_{\mathfrak R_{v}}\|x\|_{\mathfrak R_v}+\|B_\theta\|_{\mathfrak R_v}\|T_v\|_{\mathfrak R_v}\|x\|_{\mathfrak R_v}\Big)\\
&\leq& C'\|x\|_{\mathfrak R_v}.
\end{eqnarray*}

Now we can bound the diagonal operator $T_n$: given $x=(x_{n-1},\ldots,x_0)\in\mathfrak R_n$ if we set $x^1=(x_{n-1},\ldots,x_v)\in\mathfrak R_k$ and $x^2=(x_{v-1},\ldots,x_0)\in\mathfrak R_v$ then
\begin{eqnarray*}
\|T_n x\|_{\mathfrak R_n}&=&\|(T_kx^1,T_vx^2)\|_{\mathfrak R_k}=\|T_kx^1-\Omega_{v,k}(T_vx^2)\|_{\mathfrak R_k}+\|T_vx^2\|_{\mathfrak R_v}\\
&\leq&\|T_kx^1-T_k\Omega_{v,k}(x^2)\|_{\mathfrak R_k}+\|T_k\Omega_{v,k}(x^2)-\Omega_{v,k}(T_vx^2)\|_{\mathfrak R_k}+\|T_v\|\,\|x^2\|_{\mathfrak R_v}\\
&\leq& \|T_k\|\,\|x^1-\Omega_{v,k}(x^2)\|_{\mathfrak R_k}+C\|x^2\|_{\mathfrak R_v}+\|T_v\|\,\|x^2\|_{\mathfrak R_v}\\
&\leq& C\Big(\|x^1-\Omega_{v,k}(x^2)\|_{\mathfrak R_k}+\|x^2\|_{\mathfrak R_v}\Big)= C\|x\|_{\mathfrak R_n}.
\end{eqnarray*}

Conversely, if $T_n$ is bounded then the commutator is bounded since given $x\in\mathfrak R_v$ we have
\begin{eqnarray*}
\|T_k\Omega_{v,k}(x)-\Omega_{v,k}(T_vx)\|_{\mathfrak R_k}&\leq &\|T_k\Omega_{v,k}(x)-\Omega_{v,k}(T_vx)\|_{\mathfrak R_k}+\|T_vx\|_{\mathfrak R_v}\\
&=&\Big\|\Big(T_k\Omega_{v,k}(x),T_vx\Big)\Big\|_{\mathfrak R_n}\\
&=&\Big\|T_n\Big(\Omega_{v,k}(x),x\Big)\Big\|_{\mathfrak R_n}\\
&\leq& \|T_n\|\, \Big\|\Big(\Omega_{v,k}(x),x\Big)\Big\|_{\mathfrak R_n}=\|T_n\|\,\|x\|_{\mathfrak R_v}.\qedhere
\end{eqnarray*}
\end{proof}

An immediate consequence is the classical result that Domain and Range spaces are interpolation spaces (see below). An abstract formulation appears in \cite[Theorem 7.1]{racsam}, but earlier versions can be traced back to \cite{rochberg,caceso,cjmr}. We include an ad hoc proof for the sake of clarity. A space $X\subset \Sigma$ is called an interpolation space for the admissible space $\mathscr F$ if every operator $\tau: \Sigma \to \Sigma$ acting on the scale that $\mathscr F$ generates is such that $\tau[X]\subset X$. Recall that given a differential $\Omega_{k,v}$ generated from $\mathscr F$, its domain is defined as the space $\Dom \Omega_{k,v} = \{x\in \mathfrak R_k: \Omega_{k,v} x \in \mathfrak R_v \}$ endowed with the quasi-norm $\|x\|_{\Dom} =\|\Omega_{k,v} x\|_{\mathfrak R_v} +\|x\|_{\mathfrak R_k}$. The range of $\Omega_{k,v}$ is defined as $\Ran \Omega_{k,v} = \{w\in\Sigma: \exists x\in \mathfrak R_k \;: w-\Omega_{k,v} x \in
\mathfrak R_v\}$ endowed with $\|w\|_{\Ran} = \inf \{\|w - \Omega_{k,v} x \|_{\mathfrak R_v}+
\|x\|_{\mathfrak R_k}: x\in \mathfrak R_k,\;\mathrm{and} w - \Omega_{k,v} x \in \mathfrak R_v\}$. One has:

\begin{prop}\label{L_dominio}
The domain $\Dom \Omega_{1,n}$ and the range $\Ran \Omega_{n,1}$ are interpolation spaces for $\mathscr F$.
\end{prop}
\begin{proof} Let $\tau$ be an operator acting on the scale and let $T:\mathscr F\To \mathscr F$ be the induced operator. We must show that if$x\in \Dom \Omega_{1,n}$ then $\tau x\in \Dom \Omega_{1,n}$, namely, that $\Omega_{1,n}(\tau x)\in\mathfrak R_{n}$. By the Commutator Theorem \ref{Commutator} we have that, for all $y\in\mathfrak R_1$, $\Omega_{1,n}(\tau y)-T_{n}\Omega_{1,n}(y)\in \mathfrak R_n$ and$$ \Big\|\Omega_{1,n}(\tau y) - T_{n}\Omega_{1,n}(y)\Big\|_{\mathfrak R_{n}}\leq C\|y\|_{\mathfrak R_1}.$$ Thus we deduce that $ \Omega_{1,n}(\tau x)=\Omega_{1,n}(\tau x)-T_{n}\Omega_{1,n}(x) + T_{n}\Omega_{1,n}(x)\in\mathfrak R_{n}$since the last term is also in $\mathfrak R_{n}$ because $T_{n}$ is bounded on $\mathfrak R_{n}$.\medskip

The proof for range spaces is similar. Given $w \in \Ran \Omega_{n,1}$, there exist $x\in\mathfrak R_{n}$ such that $w-\Omega_{n,1}(x)\in\mathfrak R_1$. In order to show that $\tau w\in \Ran \Omega_{n,1}$ we must find $x'\in\mathfrak R_{n}$ such that $\tau w - \Omega_{n,1}(x')\in\mathfrak R_1$. Pick $x'=T_{n}(x)$ and we get
$ \tau w - \Omega_{n,1}(x')= \tau w -\Omega_{n,1}(T_{n}x) + \tau \Omega_{n,1}(x) - \tau \Omega_{n,1}(x)\in\mathfrak R_1$:
indeed, $\tau w - \tau \Omega_{n,1}(x) = \tau (w-\Omega_{n,1}(x))\in\mathfrak R_1$ by the interpolation property and $\Omega_{n,1}(T_{n}x)- \tau \Omega_{n,1}(x)\in\mathfrak R_1$ by the Commutator Theorem. \end{proof}

\end{document}